\documentclass{amsart}
\input diagrams
\usepackage{here}

\theoremstyle{plain}
 \newtheorem{thm}{Theorem}[section]
 \newtheorem{cor}{Corollary}[section]
 \newtheorem{lem}{Lemma}[section]
 \newtheorem{exmp}{Example}[section]
\theoremstyle{definition}
 \newtheorem{defn}{Definition}[section]
\theoremstyle{remark}
 \newtheorem{rem}{Remark}[section]

\newcommand{\Int}{\operatorname{Int}}

\newcommand{\Ind}{\operatorname{Ind}}

\begin{document}
\noindent
Topology Atlas Invited Contributions \textbf{8} (2003) 8 pp.
\bigskip

\title{One recent development of the $M_3$ vs.\ $M_1$ problem}
\author{Takemi Mizokami}
\address{Joetsu University of Education, Japan} 
\email{mizokami@juen.ac.jp}
\maketitle

\section{The $M_3$ vs.\ $M_1$ problem}

All spaces are assumed to be regular $T_1$.
The term CP stands for \emph{closure-preserving}.

To begin with, we give the definitions of $M_i$-spaces which were
introduced by Ceder~\cite{C} in 1960 as generalized metric spaces:

\begin{defn}
A space is an \emph{$M_1$-space} if it has a $\sigma$-CP base.
\end{defn}

\begin{defn}
A space is an \emph{$M_2$-space} if it has a $\sigma$-CP quasi-base.
\end{defn}

\begin{defn}
A space is an \emph{$M_3$-space} if it has a $\sigma$-cushioned 
pair-base.
\end{defn}

Recalling the Nagata-Smirnov metrization theorem, we easily have the 
implication:
$$ 
\text{Metric space} \to M_1 \to M_2 \to M_3 
$$
Later, Borges~\cite{B} renamed $M_3$-spaces \emph{stratifiable spaces} in 
terms of the stratification as follows:
\begin{defn}
A space $X$ is \emph{stratifiable} if there exists a function 
$$
S \colon \{\text{closed subsets of $X$}\} \times \mathbb{N} \to \tau(X),
$$ 
called the \emph{stratification} of $X$, satisfying the following:
\begin{enumerate}
\renewcommand{\labelenumi}{(\roman{enumi})}
    \item For each closed subset $H$,
    \[ \bigcap_{n}S(H,n) = \bigcap_{n}\overline{S(H,n)} = H; \]
    \item if $H \subset K$, then $S(H,n) \subset S(K,n)$ for each $n$.
\end{enumerate}
\end{defn}

As for the reverse implications, Gruenhage~\cite{G} and Junnila~\cite{J}
showed independently $M_3 \to M_2$. However, $M_3 \to M_1$? has not been
answered yet and it has become one of the most outstanding open problems
in general topology.

To this problem, many partial answers have been given, since Slaughter
first showed that any La\v{s}nev space is an $M_1$-space. In a sense, we
could say that the history of this problem is one of partial answers.
A diagram of those partial answers is given in the last part.

Since in order to state our discussion it seems better to be familiar with
Ito's method which was adopted to show that Nagata spaces or $M_3$-spaces
whose every point has a CP open neighborhood base are $M_1$, we give here
the outline of his method:

\begin{thm}[\cite{I}]
If each point of an $M_3$-space $X$ has a CP open neighborhood base, then 
every closed subset of $X$ has a CP open neighborhood base, necessarily 
$X$ is an $M_1$-space.
\end{thm}

\begin{defn}
$\mathcal{P}$ represents the class of $M_3$-spaces whose each closed subset 
has a CP open neighborhood base.
\end{defn}

The following fact due to Ito \cite{I} seems useful and applicable to our 
discussion:
If $M$ be a closed subset of an $M_3$-space, then there exists a CP 
closed neighborhood base $\mathcal{B}$ of $M$ and at the same time there 
exists a dense subset $D = \bigcup\{D_{n} \mid n \in \mathbb{N}\}$, where 
each $D_{n}$ is discrete and closed in $X$, such that
\[ B = \overline{B \cap D} \]
for each $B \in \mathcal{B}$.
By the assumption, each $p \in D_{n},\, n \in \mathbb{N}$, has an open 
neighborhood base $\mathcal{U}(p)$.
Then we can construct a CP open neighborhood base of $M$ by expanding 
each point $p \in D_{n} \cap B,\, n \in $, to members in 
$\mathcal{U}(p)|S(\{p\},n) = \{ U \cap S(\{p\},n) \mid U \in \mathcal{U}(p)\}$.

On the other hand, Tamano pointed out that the assumption on points of an 
$M_3$-spaces can be weakened as follows:
\begin{thm}[\cite{T}]
If each point $p$ of an $M_3$-space $X$, there exists a CP family $\mathcal{U}$ 
of open subsets of $X$ such that\/ $\overline{\mathcal{U}} = \{\overline{U} 
\mid U \in \mathcal{U}\}$ forms a local network at $p$ in $X$, then $X \in \mathcal{P}$.
\end{thm}

So, what kind of spaces satisfies the Tamano's criteria above?
For this question, Tamano himself gave there the following positive result:

\begin{thm}[\cite{T}]
If a space $X$ is a Baire, Fr\'{e}chet $M_3$-space, then $X$ satisfies 
the Tamano's criteria, necessarily $X \in \mathcal{P}$.
\end{thm}

In the proof, Baire property is used to induce the fact that any Baire 
$\sigma$-space has a $G_\delta$-dense metrizable subset, due to van 
Douwen.

It is quite natural to ask the following problem:
Can we delete the assumption \emph{Baire} from the above theorem, that is, 
are Fr\'{e}chet $M_3$-spaces $M_1$?

We gave the positive answer to this problem in our first paper \cite{MS1}.
To do this, we define property ($*$) as follows:
\begin{defn}[\cite{MS1}]
A space $X$ satisfies \emph{property ($*$)} if for each $p \in \partial 
O,\, O \in \tau(X)$, there exists a CP closed local network $\mathcal{B}$ at 
$p$ in $X$ such that
$$
B \subset \overline{O},
B \cap \partial O = \{p\}\ \text{and}\
\overline{B \setminus \{p\}} = B
$$ 
for each $B \in \mathcal{B}$.
\end{defn}

We note that every Fr\'{e}chet space has this property, and also that this
property generalized the Tamano's criteria in the sense of ``without the
term \emph{open sets}''.

With our criteria, we showed the following:
\begin{thm}[\cite{MS1}]
An $M_3$-space with property ($*$) belongs to $\mathcal{P}$, necessarily is $M_1$.
\end{thm}

Whether the property could work well to the $M_3$ vs.\ $M_1$ problem or
whether the property is useful to the $M_3$ vs.\ $M_1$ problem actually?
For that matter, we have to clear the essential question whether property
($*$) is what $M_3$-spaces themselves have intrinsically as their
property.

But unfortunately the answer is no, which was given by the following:

\begin{exmp}[\cite{MS1}]
There exists an $M_0$-space which does not satisfy property ($*$).
\end{exmp}

This means that property ($*$) is too strong to apply to the $M_3$ vs.\ 
$M_1$ problem; so in the second paper, we weaken it to property (P):

\begin{defn}[\cite{MS2}]
A space $X$ satisfies \emph{property (P)} if for each $p \in \partial 
O,\, O \in \tau(X)$, there exists a CP closed local network $\mathcal{B}$ at 
$p$ in $X$ such that $\overline{B \cap O} = B$ for each $B \in \mathcal{B}$.
\end{defn}

Using a similar way to the case of property ($*$), we showed the following:

\begin{thm}[\cite{MS2}]
An $M_3$-space with property (P) belongs to $\mathcal{P}$.
\end{thm}

What kind of $M_3$-spaces satisfies this property?
Do $M_3$-spaces satisfy this property?
With respect to these questions, we do not have the exact answers.
Even if a space $X \in \mathcal{P}$, we do not know whether $X$ satisfies 
property (P) or not.
Rather, this question is equivalent with the $M_3$ vs.\ $M_1$ problem 
itself, as stated below, i.e., all $M_3$-space are $M_1$ if and only if 
all $X \in \mathcal{P}$ satisfy property (P).

Since the fact that in reality we do not know whether any space in 
$\mathcal{P}$ satisfies property (P) is one defect, in a final stage we relax 
the property (P) more to property ($\delta$) as follows:

\begin{defn}[\cite{MS3}]
A space $X$ satisfies \emph{property ($\delta$)} if for each nowhere 
dense, closed subset $M$ of a space and each $p \in M$, there exists a CP 
closed local network $\mathcal{B}$ at $p$ in $X$ such that $B = \overline{B 
\setminus M}$ for each $B \in \mathcal{B}$. 
\end{defn}

Obviously, property ($*$) $\to$ property (P) $\to$ property ($\delta$).
This property is rather significant than the previous two in the sense 
that all spaces in $\mathcal{P}$ satisfy property ($\delta$) in turn.
Of course, to this case we can also show the following:

\begin{thm}[\cite{MS3}]
An $M_3$-space with property ($\delta$) belongs to $\mathcal{P}$.
\end{thm}

From it, we can characterize the class $\mathcal{P}$ in terms of this property:

\begin{cor}[\cite{MS3}]
A space $X$ belongs to $\mathcal{P}$ if and only if $X$ is an $M_3$-space with 
property ($\delta$).
\end{cor}

Here, we give lemmas needed for the proof.

\begin{lem}[\cite{D}]
An $M_3$-space is a $K_1$-space in the sense of van Douwen.
\end{lem}

\begin{lem}
Let $\mathcal{B}$ be a CP family of closed subsets of an $M_3$-space $X$.
Then there exists a pair $\langle \mathcal{F},\mathcal{V} \rangle$ of families 
satisfying the following:
\begin{enumerate}
\renewcommand{\labelenumi}{\textup{(\roman{enumi})}}
    \item $\mathcal{F}$ is a $\sigma$-discrete closed cover of $X$;
    \item $\mathcal{V} = \{V(F) \mid F \in \mathcal{F}\}$ is a point-finite, 
$\sigma$-discrete open cover of $X$ such that $F \subset V(F)$ for each 
$F \in \mathcal{F}$;
    \item for each $F \in \mathcal{F}$ and $B \in \mathcal{B}$, $F \cap B \ne 
\emptyset$ implies $F \subset B$ and $F \cap B = \emptyset$ implies $V(F) 
\cap B = \emptyset$.
\end{enumerate}
\end{lem}

\begin{lem}
Let $\mathcal{B}$ be a CP family of closed subsets of an $M_3$-space $X$.
Then there exist families $\mathcal{B}(B),\, B \in \mathcal{B}$, of subsets of $X$ 
satisfying the following:
\begin{enumerate}
\renewcommand{\labelenumi}{\textup{(\roman{enumi})}}
    \item For each $B \in \mathcal{B}$, $\mathcal{B}(B)$ is a closed neighborhood 
base of $B$ in $X$;
    \item $\bigcup\{\mathcal{B}(B) \mid B \in \mathcal{B}\}$ is CP in $X$.
\end{enumerate}
\end{lem}

\begin{lem}
Let $M$ be a closed subset of an $M_3$-space $X$ and let $\mathcal{B}$ be a CP 
family of closed subsets of $X$.
Then there exist families $\mathcal{W}(B),\, B \in \mathcal{B}$, of subsets of $X$ 
satisfying the following:
\begin{enumerate}
\renewcommand{\labelenumi}{\textup{(\roman{enumi})}}
    \item $\bigcup\{\mathcal{W}(B) \mid B \in \mathcal{B}\}$ is a CP family of closed 
subsets of $X$;
    \item for each $B \in \mathcal{B}$, 
$\mathcal{W}(B)|M = \{W \cap M \mid W \in \mathcal{W}(B)\}$ and 
$\mathcal{W}(B)|(X \setminus M) =  \{W \setminus M \mid W \in \mathcal{W}(B)\}$ 
is a closed neighborhood base of $B \setminus M$ in $X \setminus M$. 
\end{enumerate}
\end{lem}

\begin{lem}
Let $M$ be a closed subset of an $M_3$-space $X$ and let $\mathcal{B}$ be a CP 
family of closed subsets of $X$.
Then there exists a family $\{S(B) \mid B \in \mathcal{B}\}$ of open subsets of 
$X \setminus M$ satisfying the following:
\begin{enumerate}
\renewcommand{\labelenumi}{\textup{(\roman{enumi})}}
    \item For each $B \in \mathcal{B}$, $B \setminus M \subset S(B)$;
    \item for any $\mathcal{B}' \subset \mathcal{B}$,
    \[ \overline{\bigcup\{S(B) \mid B \in \mathcal{B}'\}} \cap M \subset 
\left(\bigcup\mathcal{B}'\right) \cap M. \]
\end{enumerate}
\end{lem}

\begin{lem}
Let $M$ be a closed, nowhere dense subset of an $M_3$-space $X$ with 
property ($\delta$) and let $\mathcal{B}$ be a CP family of closed subsets of $X$.
Then there exist families $\{\mathcal{W}(B) \mid B \in \mathcal{B}\}$ of subsets of 
$X$ satisfying the following:
\begin{enumerate}
\renewcommand{\labelenumi}{\textup{(\roman{enumi})}}
    \item $\mathcal{W} = \bigcup\{\mathcal{W}(B) \mid B \in \mathcal{B}\}$ is a CP family 
of closed subsets of $X$;
    \item if $B \subset O$, where $B \in \mathcal{B},\, O \in \tau(X)$, then 
there exists $W \in \mathcal{W}(B)$ such that $B \subset W \subset O$;
    \item for each $W \in \mathcal{W}$,
    \[ \overline{\Int(W \setminus M)} \cap M = W \cap M. \]
\end{enumerate}
\end{lem}

The next corollary shows how the properties (P) and ($\delta$) are 
related to the $M_3$ vs.\ $M_1$ problems:

\begin{cor}[\cite{MS3}]
TFAE:
\begin{enumerate}
\renewcommand{\labelenumi}{\textup{(\roman{enumi})}}
    \item All $M_3$-spaces are $M_1$.
    \item All $M_3$-spaces belong to $\mathcal{P}$.
    \item All spaces in $\mathcal{P}$ satisfy property (P).
    \item All $M_3$-spaces satisfy property ($\delta$).
\end{enumerate}
\end{cor}

We consider what kind of spaces satisfy property ($\delta$).
To settle one possibility, we introduce the notion of $\delta$-order, 
which is one variation of sequential order.

\begin{defn}[\cite{MS3}]
Let $A$ be a subset of a space $X$.
We introduce a operator $[\cdot]$ as follows:
\[ [A] = \left\{p \in X \mid \mbox{there exists}\ B \subset A\ 
\mbox{such that}\ \overline{B} = B \cup \{p\}\right\}. 
\]
Let $A_{0} = A$ and suppose $A_{\beta},\, \beta < \alpha$, are defined.
Then we define $A_{\alpha}$ as follows:\\
If $\alpha$ is a limit ordinal, then
\[ A_{\alpha} = \bigcup\{A_{\beta} \mid \beta < \alpha\}, \]
and if $\alpha$ is isolated, then
\[ A_{\alpha} = [A_{\alpha-1}]. \]
Define $\delta(X)$ to be the least $\alpha$ such that $\overline{A} = 
A_{\alpha}$ for any $A \subset X$ and call it the \emph{$\delta$-order 
of} $X$.
\end{defn}

We remark two points.

\begin{rem}
If a space $X$ satisfies property ($*$), then $\overline{O} = [O]$ for 
any $O \in \tau(X)$.
\end{rem}

\begin{rem}
If a space has the sequential order, then it has the $\delta$-order.
But the converse is not true.
\end{rem}

We can show that every sequential $M_3$-space has property ($\delta$).
More strictly, we have the following:

\begin{thm}[\cite{MS3}]
If a space $X$ is an $M_3$-space with the $\delta$-order $\delta(X)$, 
then $X$ satisfies property ($\delta$); consequently $X \in \mathcal{P}$.
\end{thm}

\begin{cor}
If $X$ is a k-$M_3$-space, then $X \in \mathcal{P}$.
\end{cor}

At this stage, this is the strongest result to the $M_3$ vs.\ $M_1$ 
problem among the known results.

To characterize spaces with the $\delta$-order, we introduce the 
following definition:

\begin{defn}[\cite{MS3}]
A space $X$ is a \emph{$\delta$-space} when a subset $F$ is closed in 
$X$ if and only if $F$ has the property that $B \subset F$ and 
$\overline{B} = B \cup \{p\}$ imply $p \in F$.
A space $L$ is an \emph{almost discrete space} if all points of $L$ 
except one point are isolated in $L$.
\end{defn}
(After we have announced the definition of $\delta$-spaces, we are 
informed by the referee of this article that this property was originally 
introduced by A. Pultr and A. Tozzi under the name of spaces with the 
property of Weakly Approximation by Points (= WAP spaces) 
\cite{MR94f:18003}.)

\begin{thm}[\cite{MS3}]
TFAE:
\begin{enumerate}
\renewcommand{\labelenumi}{\textup{(\roman{enumi})}}
    \item $X$ is a $\delta$-space.
    \item $X$ has the $\delta$-order.
    \item $X$ is the image of $\oplus\{L_{\lambda} \mid \lambda \in 
\Lambda\}$ with each $L_{\lambda}$ almost discrete under a quotient 
mapping $\varphi$ such that $\varphi(L_{\lambda})$ is a closed 
$\delta$-space for each $\lambda \in \Lambda$.
\end{enumerate}
\end{thm}

\begin{thm}[\cite{MS3}]
For a pointwise perfect space $X$, TFAE:
\begin{enumerate}
\renewcommand{\labelenumi}{\textup{(\roman{enumi})}}
    \item $X$ is a $\delta$-space.
    \item $X$ has the $\delta$-order.
    \item $X$ is the image of $\oplus\{L_{\lambda} \mid \lambda \in 
\Lambda\}$ with each $L_{\lambda}$ an almost discrete, $M_{0}$-space 
under a quotient mapping $\varphi$ such that $\varphi(L_{\lambda})$ is a 
closed $\delta$-space for each $\lambda \in \Lambda$.
\end{enumerate}
\end{thm}

\newpage
\section{The diagrams}

\begin{figure}[H]\label{fig1}
\begin{diagram}[small,tight]
				&	&&				&	&				&	&&	&\text{Baire, Fr\'echet $M_3$-space}&&	&	\\
				&	&&				&	&\text{La\v{s}nev space}	&	&&	&\dTo				&&	&\text{Nagata space}\\
				&	&&				&	&\dTo				&\rdTo(4,3)&&	&				&&\ldTo(3,3)&	\\
\text{$M_3$, $\sigma$-discrete space}&	&&				&	&\text{L-space}			&	&&	&				&&	&	\\
\dTo				&	&&				&	&\dTo				&	&&	&\text{Fr\'echet $M_3$-space}	&&	&	\\
\text{$M_3$, $F_\sigma$-metrizable space}&&&				&	&\text{free L-space}		&	&&	&\dTo				&&	&	\\
				&\rdTo(3,2)&&				&\ldTo	&				&	&&	&				&&	&	\\
				&	&&\text{$M_3$, $\mu$-space}	&	&				&	&&	&\text{sequential $M_3$-space}	&&	&	\\
				&	&&\dTo				&	&				&	&&\ldTo(4,4)&				&&	&	\\
				&	&&\text{perfect image of an $M_0$-space}&&				&	&&	&				&&	&	\\
				&	&&				&\rdTo	&				&	&&	&				&&	&	\\
				&	&&				&	&\text{$M_3$-space with property (P)}&	&&	&				&&	&	\\
				&	&&				&	&\dTo				&	&&	&				&&	&	\\
				&	&&				&	&\text{$M_1$-space}		&	&&	&				&&	&	\\
\end{diagram}
\smallskip
\caption{}
\end{figure}

\newarrow{Corresponds}{<}{-}{-}{-}{>}
\begin{figure}[H]\label{fig2}
\begin{diagram}[small]
\text{Nagata space}			&	&				\\
\dTo					&	&				\\
\text{sequential $M_3$-space}		&\rTo	&\text{$M_3$-space with property (P)}\\
\dTo					&	&\dTo				\\
\text{$M_3$-space with the $\delta$-order}&\rTo	&\text{$M_3$-space with property ($\delta$)}\\
\dCorresponds				&	&\dCorresponds			\\
\text{$M_3$-$\delta$-space}		&\rTo	&\mathcal{P}			\\
\end{diagram}
\smallskip
\caption{}
\end{figure}

\newpage
\section{Related problems}

We list here other open problems that are equivalent or related to the 
$M_3$ vs.\ $M_1$ problem.
They are already proposed elsewhere by many researchers.
\begin{enumerate}
\renewcommand{\labelenumi}{(P\theenumi)}
\item 
Is any $M_{3}$-space $M_{1}$?
\item 
Is any (closed) subspace of an $M_{1}$-space $M_{1}$?
\item 
Is any closed (perfect) image of an $M_1$-space $M_1$?
\item 
Does any point of an $M_{1}$-space have a CP open neighborhood base?
\item 
Is any adjunction space $X\cup_{f}Y$ for $M_{1}$-spaces $X$ and $Y$ $M_{1}$?
\item 
Is any $M_{1}$-space the perfect image of an $M_{1}$-space with $\Ind = 0$?
\item 
Is any $M_3$-space the perfect image of an $M_3$-space with $\Ind = 0$?
\item 
Does $\mathcal{E}\mathcal{M}_3 \subset \{\mbox{M}_1\mbox{-spaces}\}$?
\item 
Does $\{\mbox{M}_1\mbox{-spaces}\} \subset \mathcal{E}\mathcal{M}_3$?
\item 
For any $M_1$-space $X$, can we characterize $\dim X \le n$ or $\Ind X \le 
n$ by the fact that there exists a $\sigma$-CP base $\mathcal{W}$ for X such 
that $\dim \partial W \le n-1$ or $\Ind \partial W \le n-1$, 
respectively, for any $W \in {\mathcal{W}}$?
\item 
If an $M_1$-space $X$ has  $\Ind = 0$, then is $X$ an $M_0$-space?
\item 
Does any $M_3$-space have an M-structure?
\item 
Does there exist a subclass $\mathcal{C}$ of $M_1$-spaces satisfying the 
following topological operations:
\begin{enumerate}
\renewcommand{\labelenumii}{(\roman{enumii})}
\item $\mathcal{C}$ is hereditary;
\item $\mathcal{C}$ is countably productive;
\item $\mathcal{C}$ is preserved by closed mappings?
\end{enumerate}
\end{enumerate}

\begin{figure}[H]
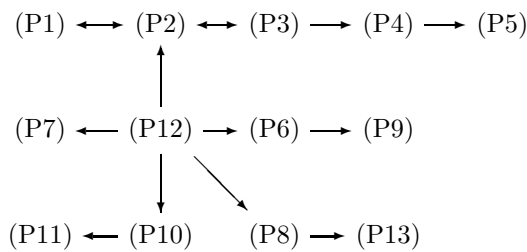
\label{fig3}
\begin{diagram}[small]
\text{(P1)}	&\rCorresponds	&\text{(P2)}	&\rCorresponds	&\text{(P3)}	&\rTo	&\text{(P4)}	&\rTo	&\text{(P5)}\\
		&		&\uTo		&		&		&	&		&	&	\\
\text{(P7)}	&\lTo		&\text{(P12)}	&\rTo		&\text{(P6)}	&\rTo	&\text{(P9)}	&	&	\\
		&		&\dTo		&\rdTo		&		&	&		&	&	\\
\text{(P11)}	&\lTo		&\text{(P10)}	&		&\text{(P8)}	&\rTo	&\text{(P13)}	&	&	\\
\end{diagram}
\smallskip
\caption{Relationships among the problems (P$n$)}
\end{figure}

Finally, we give a diagram (Figure \ref{fig3}) among all the problems stated here, where A 
$\longrightarrow$ B means that if A is positively solved, then so is B.

\begin{thebibliography}{10}

\bibitem{B}
C.~J.~R. Borges, \emph{On stratifiable spaces}, Pacific J. Math. \textbf{17}
  (1966), 1--16. \MR{32 \#6409}

\bibitem{C}
J.~G. Ceder, \emph{Some generalizations of metric spaces}, Pacific J. Math.
  \textbf{11} (1961), 105--125. \MR{24 \#A1707}

\bibitem{D}
E.~K.~van Douwen, \emph{Simultaneous extension of continuous functions}, Ph.D.
  thesis, Vrije Univesiteit, Amsterdam, 1975, p.~75 pp.

\bibitem{G}
G.~Gruenhage, \emph{Stratifiable spaces are ${M}\sb{2}$}, Topology Proc.
  \textbf{1} (1976), 221--226. \MR{56 \#6614}

\bibitem{I}
M.~It{\=o}, \emph{${M}_3$-spaces whose every point has a closure preserving
  outer base are ${M}_1$}, Topology Appl. \textbf{19} (1985), no.~1, 65--69.
  \MR{86g:54039}

\bibitem{J}
H.~J.~K. Junnila, \emph{Neighbornets}, Pacific J. Math. \textbf{76} (1978),
  no.~1, 83--108. \MR{58 \#2734}

\bibitem{MS1}
T.~Mizokami and N.~Shimane, \emph{On {F}r\'echet ${M}_3$-spaces}, Math. Japon.
  \textbf{50} (1999), no.~3, 391--399. \MR{2000j:54032}

\bibitem{MS2}
\bysame, \emph{On the ${M}_3$ versus ${M}_1$ problem}, Topology Appl.
  \textbf{105} (2000), no.~1, 1--13. \MR{2001f:54032}

\bibitem{MS3}
T.~Mizokami, N.~Shimane, and Y.~Kitamura, \emph{A characterization of a certain
  subclass of ${M}_1$-spaces}, JP J. Geom. Topol. \textbf{1} (2001), no.~1,
  37--51. \MR{1 876 154}

\bibitem{MR94f:18003}
A.~Pultr and A.~Tozzi, \emph{Equationally closed subframes and representation
  of quotient spaces}, Cahiers Topologie G\'eom. Diff\'erentielle Cat\'eg.
  \textbf{34} (1993), no.~3, 167--183. \MR{94f:18003}

\bibitem{T}
K.~Tamano, \emph{$\mu$-spaces, stratifiable spaces and mosaical collections},
  Math. Japon. \textbf{34} (1989), no.~3, 483--496. \MR{90e:54062}

\end{thebibliography}
\providecommand{\bysame}{\leavevmode\hbox to3em{\hrulefill}\thinspace}
\providecommand{\MR}{\relax\ifhmode\unskip\space\fi MR }
\providecommand{\MRhref}[2]{%
  \href{http://www.ams.org/mathscinet-getitem?mr=#1}{#2}
}
\providecommand{\href}[2]{#2}

\end{document}